\documentclass[graybox,envcountsect]{svmult}

\usepackage{mathptmx}       
\usepackage{helvet}         
\usepackage{courier}        
\usepackage{type1cm}        
\usepackage{makeidx}         
\usepackage{graphicx}        
\usepackage{multicol}        
\usepackage[bottom]{footmisc}

\usepackage{amsmath,amssymb,latexsym,amsbsy,natbib}

\makeindex             

\usepackage{amsfonts, amsbsy, calrsfs}
\usepackage{natbib}

\makeatletter

\newcommand{\Rmnum}[1]{\expandafter\@slowromancap\romannumeral #1@}
\makeatother

\renewcommand{\le}{\leqslant}
\renewcommand{\ge}{\geqslant}

\DeclareMathOperator{\expec}{\mathrm{E}}

\DeclareMathOperator{\proba}{\mathrm{P}}

\newcommand{\mv}{\boldsymbol}
\newcommand{\defn}{\emph}
\newcommand{\RR}{\mathbb{R}}
\newcommand{\simplex}{\Delta_{d-1}}
\newcommand{\diff}{\mathrm{d}}

\bibpunct{[}{]}{,}{,}{,}{,}

\begin{document}

\title*{Extreme-Value Copulas}
\author{Gordon Gudendorf and Johan Segers}
\institute{Gordon Gudendorf \at Institut de statistique, Universit\'e catholique de Louvain, Louvain-la-Neuve, Belgium; \email{gordon.gudendorf@uclouvain.be}
\and Johan Segers \at  Institut de statistique, Universit\'e catholique de Louvain, Louvain-la-Neuve, Belgium; \email{johan.segers@uclouvain.be}}
%
%


\maketitle

\abstract{Being the limits of copulas of componentwise maxima in independent random samples, extreme-value copulas can be considered to provide appropriate models for the dependence structure between rare events. Extreme-value copulas not only arise naturally in the domain of extreme-value theory, they can also be a convenient choice to model general positive dependence structures. The aim of this survey is to present the reader with the state-of-the-art in dependence modeling via extreme-value copulas. Both probabilistic and statistical issues are reviewed, in a nonparametric as well as a parametric context.} 
\bigskip

\noindent This version: December 03, 2009

\section{Introduction}

In various domains, as for example finance, insurance or environmental science, joint extreme events can have a serious impact and therefore need careful modeling. Think for instance of daily water levels at two different locations in a lake during a year. Calculation of the probability that there is a flood exceeding a certain benchmark requires knowledge of the joint distribution of maximal heights during the forecasting period. This is a typical field of application for extreme-value theory. In such situations, extreme-value copulas can be considered to provide appropriate models for the dependence structure between exceptional events. 

One of the first applications of bivariate extreme-value analysis must be due to Gumbel and Goldstein \cite{GG64}, who analyze the maximal annual discharges of the Ocmulgee River in Georgia at two different stations, a dataset that has been taken up again in \cite{RVF01}. The joint behavior of extreme returns in the foreign exchange rate market is investigated in \cite{S99}, whereas the comovement of equity markets characterized by high volatility levels is studied in \cite{LS01}. An application in the insurance domain can be found in \cite{CDL03}.

Extreme-value copulas not only arise naturally in the domain of extreme events, but they can also be a convenient choice to model data with positive dependence. An advantage with respect to the much more popular class of Archimedean copulas, for instance, is that they are not symmetric. Incidentally, a hybrid class containing both the Archimedean and the extreme-value copulas as a special case are the Archimax copulas~\cite{CFG00}.

The aim of this survey is to present the reader with the state-of-the-art in dependence modeling via extreme-value copulas. Definition, origin, and basic properties of extreme-value copulas are presented in Section~\ref{S:foundations}. A number of useful and popular parametric families are reviewed in Section~\ref{S:parametric}. Section~\ref{S:coefficients} provides a discussion of the most important dependence coefficients associated to extreme-value copulas. An overview of parametric and nonparametric inference methods for extreme-value copulas is given in Section~\ref{S:estim}. Finally, some further topics and pointers to the literature are gathered in Section~\ref{S:further}.

\section{Foundations}
\label{S:foundations}

Let $\mv{X}_i = (X_{i1}, \ldots, X_{id})$, $i \in \{1, \ldots, n\}$, be a sample of independent and identically distributed (iid) random vectors with common distribution function $F$, margins $F_1, \ldots, F_d$, and copula $C_F$. For convenience, assume $F$ is continuous. Consider the vector of componentwise maxima:
\begin{equation}
\label{E:max}
  \mv{M}_{n} = (M_{n,1}, \ldots, M_{n,d}), \qquad \text{where } M_{n,j} = \bigvee_{i=1}^n X_{ij},
\end{equation}
with `$\vee$' denoting maximum. Since the joint and marginal distribution functions of $\mv{M}_n$ are given by $F^n$ and $F_1^n, \ldots, F_d^n$ respectively, it follows that the copula, $C_n$, of $\mv{M}_n$ is given by
\[
  C_{n}(u_1, \ldots, u_d) = C_{F}(u_1^{1/n}, \ldots, u_d^{1/n})^n, \qquad (u_1, \ldots, u_d) \in [0,1]^d.
\]
The family of extreme-value copulas arises as the limits of these copulas $C_n$ as the sample size $n$ tends to infinity. 

\begin{definition}
\label{extv:cop}
A copula $C$ is called an \defn{extreme-value copula} if there exists a copula $C_F$ such that
\begin{equation}
\label{E:DA}
  C_{F}(u_1^{1/n}, \ldots, u_d^{1/n})^n \to C(u_1, \ldots, u_d) \qquad (n \to \infty)
\end{equation}
for all $(u_1, \ldots, u_d) \in [0,1]^{d}$. The copula $C_F$ is said to be in the \defn{domain of attraction} of $C$.
\end{definition}

Historically, this construction dates back at least to \cite{Deheuvels84, Galambos78}. The representation of extreme-value copulas can be simplified using the concept of max-stability.

\begin{definition} 
\label{D:maxstable}
A $d$-variate copula $C$ is \defn{max-stable} if it satisfies the relationship
\begin{equation}
\label{E:maxstable}
  C(u_1, \ldots, u_d) = C(u_1^{1/m}, \ldots, u_d^{1/m})^m
\end{equation}
for every integer $m \ge 1$ and all $(u_1, \ldots, u_d) \in [0,1]^{d}$.
\end{definition}

From the previous definitions, it is trivial to see that a max-stable copula is in its own domain of attraction and thus must be itself an extreme-value copula. The converse is true as well.

\begin{theorem}
\label{T:EVC:maxstable}
A copula is an extreme-value copula if and only if it is max-stable.
\end{theorem}

The proof of Theorem~\ref{T:EVC:maxstable} is standard: for fixed integer $m \ge 1$ and for $n = mk$, write
\[
  C_n(u_1, \ldots, u_d) = C_k(u_1^{1/m}, \ldots, u_d^{1/m})^m.
\]
Let $k$ tend to infinity on both sides of the previous display to get \eqref{E:maxstable}.

\begin{figure}[h]
\begin{center}
\begin{tabular}{cc}
\includegraphics[width=0.45\textwidth]{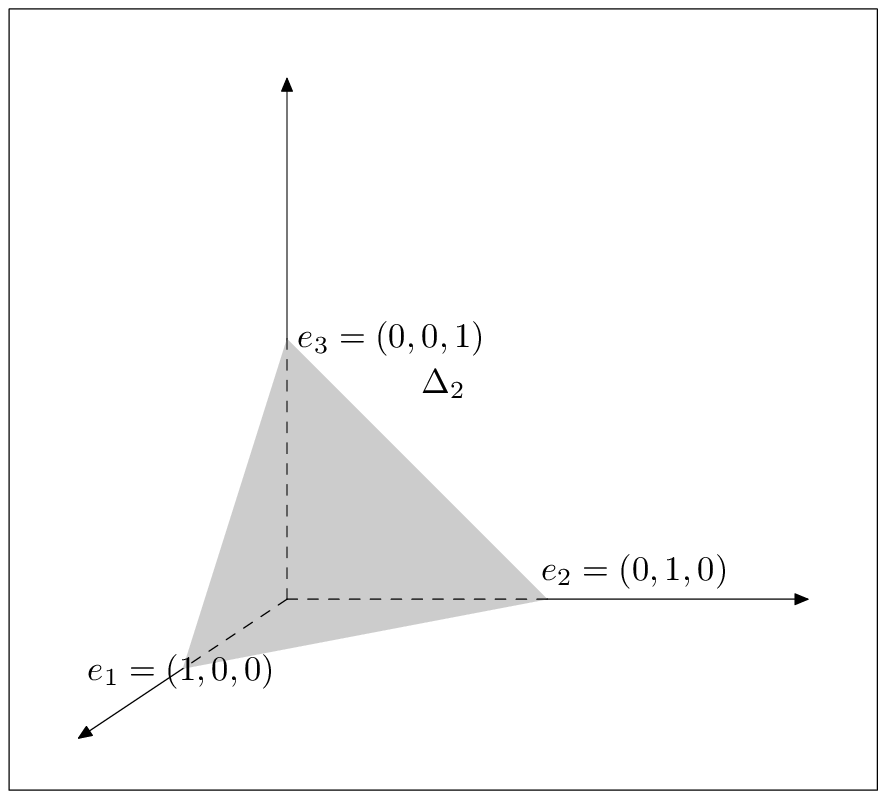}& 
\includegraphics[width=0.45\textwidth]{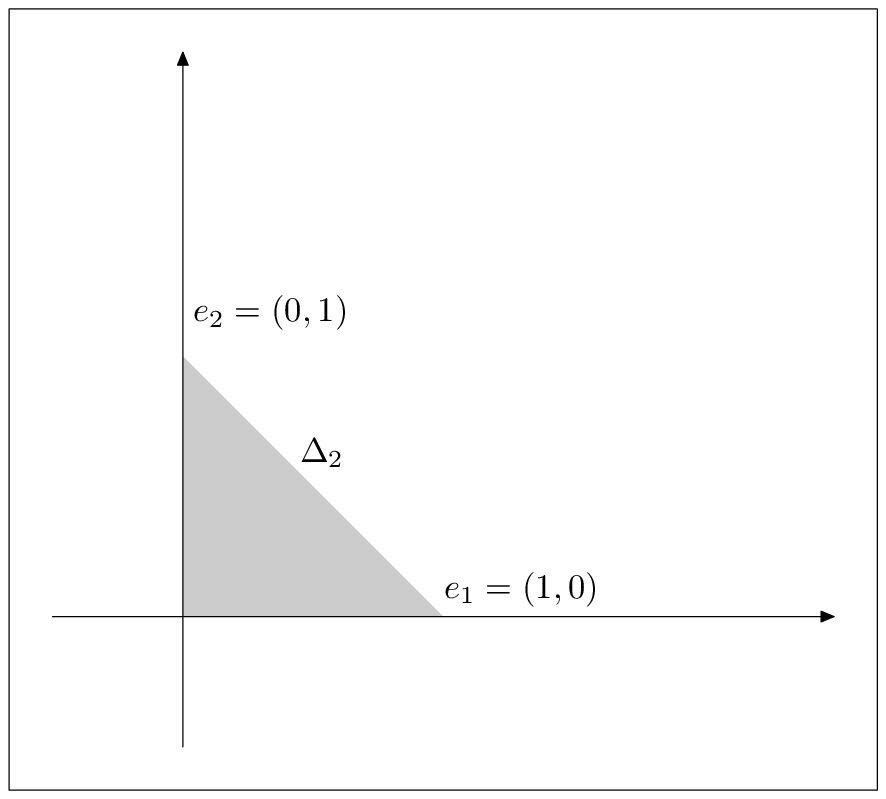}\\
\end{tabular}
\end{center}
\caption{\label{F:simplex} On the left side $\Delta_{2}$ is represented in $\RR^{3}$, which is equivalent to the representation in $\RR^{2}$ on the right side.}
\end{figure}

By definition, the family of extreme-value copulas coincides with the set of copulas of \defn{extreme-value distributions}, that is, the class of limit distributions with non-degenerate margins of
\[
  \biggl( \frac{M_{n,1} - b_{n,1}}{a_{n,1}}, \ldots, \frac{M_{n,d} - b_{n,d}}{a_{n,d}} \biggr)
\]
with $M_{n,j}$ as in \eqref{E:max}, centering constants $b_{n,j} \in \RR$ and scaling constants $a_{n,j} > 0$. Representations of extreme-value distributions then yield representations of extreme-value copulas. Let $\simplex = \{ (w_1, \ldots, w_d) \in [0, \infty)^d : \sum_j w_j = 1 \}$ be the unit simplex in $\RR^d$; see Figure~\ref{F:simplex}. The following theorem is adapted from \cite{Pickands81}, which is based in turn on \cite{HR77}. 

\begin{theorem}
\label{T:EVC:H}
A $d$-variate copula $C$ is an extreme-value copula if and only if there exists a finite Borel measure $H$ on $\simplex$, called \defn{spectral measure}, such that
\[
  C(u_1, \ldots, u_d) = \exp \bigl( - \ell( - \log u_1, \ldots, - \log u_d) \bigr), \qquad (u_1, \ldots, u_d) \in (0, 1]^d,
\]
where the \defn{tail dependence function} $\ell : [0, \infty)^d \to [0, \infty)$ is given by
\begin{equation}
\label{E:ell}
  \ell(x_1, \ldots, x_d) = \int_{\simplex} \bigvee_{j=1}^d (w_j x_j) \; dH(w_1, \ldots, w_d), \qquad (x_1, \ldots, x_d) \in [0, \infty)^d.
\end{equation}
The spectral measure $H$ is arbitrary except for the $d$ moment constraints
\begin{equation}
\label{E:constraints}
  \int_{\simplex} w_j \; dH(w_1, \ldots, w_d) = 1, \qquad j \in \{1, \ldots, d\}.
\end{equation}
\end{theorem}

The $d$ moment constraints on $H$ in \eqref{E:constraints} stem from the requirement that the margins of $C$ be standard uniform. They imply that $H(\simplex) = d$.

By a linear expansion of the logarithm and the exponential function, the domain-of-attraction equation~\eqref{E:DA} is equivalent to
\begin{align}
\label{E:DA:ell}
  \lim_{t \downarrow 0} t^{-1} \bigl( 1 - C_F(1 - tx_1, \ldots, 1 - tx_d) \bigr) 
  &= - \log C ( e^{-x_1}, \ldots, e^{-x_d} ) \nonumber \\
  &= \ell(x_1, \ldots, x_d)
\end{align}
for all $(x_1, \ldots, x_d) \in [0, \infty)^d$; see for instance \cite{DH98}. The tail dependence function $\ell$ in \eqref{E:ell} is convex, homogeneous of order one, that is $\ell(cx_1, \ldots, cx_d) = c \, \ell(x_1, \ldots, x_d)$ for $c > 0$, and satisfies $\max(x_1, \ldots, x_d) \le \ell(x_1, \ldots, x_d) \le x_1 + \cdots + x_d$ for all $(x_1, \ldots, x_d) \in [0, \infty)^d$. By homogeneity, it is characterized by the \defn{Pickands dependence function} $A : \simplex \to [1/d, 1]$, which is simply the restriction of $\ell$ to the unit simplex:
\[
  \ell(x_1, \ldots, x_d) = (x_1 + \cdots + x_d) \, A(w_1, \ldots, w_d) \qquad \text{where } \qquad w_j = \frac{x_j}{x_1 + \cdots + x_d},
\]
for $(x_1, \ldots, x_d) \in [0, \infty)^d \setminus \{0\}$. The extreme-value copula $C$ can be expressed in terms of $A$ via
\[
  C(u_1, \ldots, u_d)
  = \exp \left\{ \left( \sum_{j=1}^{d} \log u_{j} \right) 
    A \left(  \frac{\log u_{1}}{ \sum_{j=1}^{d} \log u_{j} }, \ldots ,  \frac{\log u_{d}}{ \sum_{j=1}^{d} \log u_{j} }  \right) \right\}.
\]
The function $A$ is convex as well and satisfies $\max(w_1, \ldots, w_d) \le A(w_1, \ldots, w_d) \le 1$ for all $(w_1, \ldots, w_d) \in \simplex$. However, these properties do not characterize the class of Pickands dependence functions unless $d = 2$, see for instance the counterexample on p.~257 in \cite{BGTS04}. 

In the bivariate case, we identify the unit simplex $\Delta_1 = \{ (1-t, t) : t \in [0, 1] \}$ in $\RR^2$ with the interval $[0, 1]$.

\begin{theorem}
\label{T:EVC:A}
A bivariate copula $C$ is an extreme-value copula if and only if
\begin{equation}
\label{E:EVC:A}
  C(u, v) = (uv)^{A( \log(v) / \log(uv) )}, \qquad (u, v) \in (0, 1]^2 \setminus \{(1, 1)\},
\end{equation}
where $A : [0, 1] \to [1/2, 1]$ is convex and satisfies $t \vee (1-t) \le A(t) \le 1$ for all $t \in [0, 1]$.
\end{theorem}

\begin{figure}[h]
\sidecaption
\includegraphics[width=0.45\textwidth]{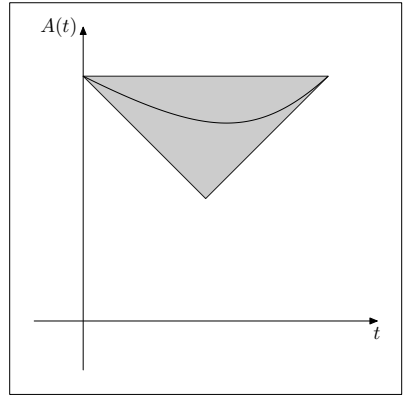} 
\caption{A typical Pickands dependence function $A$ together with the region $t \vee (1-t) \le A(t) \le 1$ in Theorem \ref{T:EVC:A}.}
\label{shadreg}
\end{figure}

It is worth stressing that in the bivariate case, any function $A$ satisfying the two constraints from Theorem~\ref{T:EVC:A} corresponds to an extreme-value copula. These functions lie in the shaded area of Figure~\ref{shadreg}; in particular, $A(0) = A(1) = 1$. 
  
The upper and lower bounds for $A$ have special meanings: the upper bound $A(t) = 1$ corresponds to independence, $C(u, v) = uv$, whereas the lower bound $A(t) = t \vee (1-t)$ corresponds to perfect dependence (comonotonicity) $C(u, v) = u \wedge v$. In general, the inequality $A(t) \le 1$ implies $C(u, v) \ge uv$, that is, extreme-value copulas are necessarily \defn{positive quadrant dependent}.

\section{Parametric models}
\label{S:parametric}

%


By Theorems~\ref{T:EVC:H} and \ref{T:EVC:A}, the class of extreme-value copulas is infinite-dimensional. Parametric submodels can be constructed in a number of ways: by calculating the limit $\ell$ in \eqref{E:DA:ell} for a given initial copula $C_F$; by specifying a spectral measure $H$; in dimension $d = 2$, by constructing a Pickands dependence function $A$. In this section, we employ the first of these methods to introduce some of the more popular families. For more extensive overviews, see e.g.\ \cite{BGTS04, NK00}.

\subsection{Logistic model or Gumbel--Hougaard copula}

\newcommand{\inv}{^\leftarrow}
Consider the Archimedean copula
\begin{equation}
\label{E:Arch}
  C_\phi(u_1, \ldots, u_d) = \phi\inv \bigl( \phi(u_1) + \cdots + \phi(u_d) \bigr), \qquad (u_1, \ldots, u_d) \in [0, 1]^d
\end{equation}
with generator $\phi : [0, 1] \to [0, \infty]$ and inverse $\phi\inv(t) = \inf \{ u \in [0, 1] : \phi(u) \le t \}$; the function $\phi$ should be strictly decreasing and convex and satisfy $\phi(1) = 0$, and $\phi\inv$ should be $d$-monotone on $(0, \infty)$, see \cite{McNN09}. 

If the following limit exists,
\begin{equation}
\label{E:logistic:theta}
  \theta = - \lim_{s \downarrow 0} \frac{s \, \phi'(1-s)}{\phi(1-s)} \in [1, \infty]
\end{equation}
then the domain-of-attraction condition~\eqref{E:DA:ell} is verified for $C_F$ equal to $C_\phi$, the tail dependence function being
\begin{equation}
\label{E:logistic:ell}
  \ell(x_1, \ldots, x_d) = 
  \begin{cases}
    (x_1^\theta + \cdots + x_d^\theta)^{1/\theta} & \text{if $1 \le \theta < \infty$}, \\
    x_1 \vee \cdots \vee x_d & \text{if $\theta = \infty$,}
  \end{cases}
\end{equation}
for $(x_1, \ldots, x_d) \in [0, \infty)^d$; see \cite{CFG00, CS09}. The range $[1, \infty]$ for the parameter $\theta$ in \eqref{E:logistic:theta} is not an assumption but rather a consequence of the properties of $\phi$. The parameter $\theta$ measures the degree of dependence, ranging from independence ($\theta = 1$) to complete dependence ($\theta = \infty$).

The extreme-value copula associated to $\ell$ in \eqref{E:logistic:ell} is
\[
  C(u_1, \ldots, u_d) = \exp \bigl\{ - \bigl( (- \log u_1)^\theta + \cdots + (- \log u_d)^\theta \bigr)^{1/\theta} \},
\]
known as the \emph{Gumbel--Hougaard} or \emph{logistic} copula. Dating back to Gumbel \cite{G60, Gumbel61}, it is (one of) the oldest multivariate extreme-value models. It was discovered independently in survival analysis \cite{Cr89, Hougaard86}. It happens to be the only copula that is at the same time Archimedean and extreme-value \cite{GR89}.

The bivariate asymmetric logistic model introduced in \cite{Tawn88} adds further flexibility to the basic logistic model. Multivariate extensions of the asymmetric logistic model were studied already in \cite{McF78} and later in \cite{CT91, Joe94}. These distributions can be generated via mixtures of certain extreme-value distributions over stable distributions, a representation that yields large possibilities for modelling that have yet begun to be explored \cite{FNR09, TGNVC09}.

\subsection{Negative logistic model or Galambos copula}

Let $\hat{C}_\phi$ be the survival copula of the Archimedean copula $C_\phi$ in \eqref{E:Arch}. Specifically, if $C_\phi$ is the distribution function of the random vector $(U_1, \ldots, U_d)$, then $\hat{C}_\phi$ is the distribution function of the random vector $(1-U_1, \ldots, 1-U_d)$. If the following limit exists,
\begin{equation}
\label{E:neglogistic:theta}
  \theta = - \lim_{s \downarrow 0} \frac{\phi(s)}{s \, \phi'(s)} \in [0, \infty]
\end{equation}
then the domain-of-attraction condition~\eqref{E:DA:ell} is verified for $C_F$ equal to $\hat{C}_\phi$, the tail dependence function being
\[
  \ell(x_1, \ldots, x_d) =
  \begin{cases}
  x_1 + \cdots + x_d & \text{if $\theta = 0$}, \\
  \displaystyle x_1 + \cdots + x_d - \sum_{I \subset \{1, \ldots, d\} \atop |I| \ge 2} (-1)^{|I|} 
    \bigl( {\textstyle\sum_{i \in I} x_i^{-\theta}} \bigr)^{-1/\theta} & \text{if $0 < \theta < \infty$}, \\
  x_1 \vee \cdots \vee x_d & \text{if $\theta = \infty$},
  \end{cases}
\]
for $(x_1, \ldots, x_d) \in [0, \infty)^d$; see \cite{CFG00, CS09}. In case $0 < \theta < \infty$, the sum is over all subsets $I$ of $\{1, \ldots, d\}$ of cardinality $|I|$ at least $2$. The amount of dependence ranges from independence ($\theta = 0$) to complete dependence ($\theta = \infty$).

The resulting extreme-value copula is known as the \emph{Galambos} or \emph{negative logistic} copula, dating back to \cite{Galambos75}. Asymmetric extensions have been proposed in \cite{Joe90, Joe94}.

\subsection{H\"usler--Reiss model}

For the bivariate normal distribution with correlation coefficient $\rho$ smaller than one, it is known since \cite{S60} that the marginal maxima $M_{n,1}$ and $M_{n,2}$ are asymptotically independent, that is, the domain-of-attraction condition~\eqref{E:DA} holds with limit copula $C(u, v) = uv$. However, for $\rho$ close to one, better approximations to the copula of $M_{n,1}$ and $M_{n,2}$ arise within a somewhat different asymptotic framework. More precisely, as in \cite{HR89}, consider the situation where the correlation coefficient $\rho$ associated to the bivariate Gaussian copula $C_\rho$ is allowed to change with the sample size, $\rho = \rho_n$, in such a way that $\rho_n \to 1$ as $n \to \infty$. If
\[
  (1 - \rho_n) \log n \to \lambda^{2} \in [0, + \infty] \qquad (n \to \infty),
\]
then one can show that
\[
  C_{\rho_n}(u^{1/n}, v^{1/n})^n \to C_A(u, v) \qquad (n \to \infty), \qquad (u, v) \in [0, 1]^2,
\]
where the \emph{H\"usler--Reiss} copula $C_A$ is the bivariate extreme-value copula with Pickands dependence function
\[
  A(w) = (1-w) \, \Phi \biggl( \lambda  + \frac{1}{2 \lambda} \log \frac{1-w}{w} \biggr) 
  + w \, \Phi \biggl( \lambda  + \frac{1}{2 \lambda} \log \frac{w}{1-w} \biggr)
\]
for $w \in [0, 1]$, with $\Phi$ representing the standard normal cumulative distribution function. The parameter $\lambda$ measures the degree of dependence, going from independence ($\lambda = \infty$) to complete dependence ($\lambda = 0$).

\subsection{The t-EV copula}

In financial applications, the $t$-copula is sometimes preferred over the Gaussian copula because of the larger weight it assigns to the tails. The bivariate $t$-copula with $\nu > 0$ degrees of freedom and correlation parameter $\rho \in (-1, 1)$ is the copula of the bivariate $t$-distribution with the same parameters and is given by
\[
  C_{\nu,\rho}(u,v) 
  = \int_{-\infty}^{t^{-1}_{\nu}( u )} \int_{- \infty}^{t^{-1}_{\nu}( v )} 
    \frac{1}{\pi \nu |P|^{1/2}} \frac{\Gamma \left(  \frac{\nu}{2} + 1 \right)}{\Gamma \left( \frac{\nu }{ 2 } \right)}
    \biggl(  1 + \frac{\mv{x}^{'} P^{-1} \mv{x}}{\nu} \biggr)^{-\nu / 2 +1} \, d\mv{x},
\]
where $t_{\nu}$ represents the distribution function of the univariate $t$-distribution with $\nu$ degrees of freedom and $P$ represents the $2 \times 2$ correlation matrix with off-diagonal element $\rho$. In \cite{DM05}, it is shown that $C_{\nu,\rho}$ is in the domain of attraction of the bivariate extreme-value copula $C_A$ with Pickands dependence function
\begin{multline}
  A(w) = w \, t_{\nu + 1} (z_w) + (1-w) \, t_{\nu + 1} (z_{1-w}), \\
  \text{where } z_{w} = (1 + \nu)^{1/2} [\{w/(1-w)\}^{1/\nu} - \rho] (1 - \rho^2)^{-1/2}, \qquad w \in [0, 1].
\end{multline}
This extreme-value copula was coined the \emph{t-EV} copula. Building upon results in~\cite{AFG05, H05}, exactly the same extreme-value attractor is found in \cite{AJ07} for the more general class of (meta-)elliptical distributions whose generator has a regularly varying tail.


\section{Dependence coefficients}
\label{S:coefficients}

Let $(U, V)$ be a bivariate random vector with distribution function $C$, a bivariate extreme-value copula with Pickands dependence function $A$ as in \eqref{E:EVC:A}. As mentioned already, the inequality $A \le 1$ implies that $C(u, v) \ge uv$ for all $(u, v) \in [0, 1]^2$, that is, $C$ is positive quadrant dependent. In fact, in \cite{G00} it was shown that extreme-value copulas are \defn{monotone regression dependent}, that is, the conditional distribution of $U$ given $V = v$ is stochastically increasing in $v$ and \textit{vice versa}; see also Theorem~5.2.10 in \cite{Resnick87}. 

In particular, all measures of dependence of $C$ such as Kendall's $\tau$ or Spearman's $\rho_S$ must be nonnegative. The latter two can be expressed in terms of $A$ via
\begin{align*}
  \tau &= 4 \iint_{[0, 1]^2} C(u,v) \, \diff C(u,v) - 1 = \int_{0}^{1} \frac{t(1-t)}{A(t)} \, \diff A'(t), \\
  \rho_S &= 12 \iint_{[0, 1]^2} u v \, \diff C(u,v) - 3 = 12 \int_{0}^{1} \frac{1}{(1+A(t))^{2}} \, \diff t - 3.
\end{align*}
The Stieltjes integrator $\diff A'(t)$ is well-defined since $A$ is a convex function on $[0,1]$; if the dependence 
function $A$ is twice differentiable, it can be replaced by $A''(t) \, \diff t$. For a proof of the identities above, see for instance \cite{H03}, where it is shown that $\tau$ and $\rho_S$ satisfy $-1 + \sqrt{1 + 3 \tau } \le \rho_{S} \le \min \left( \frac{3}{2} \tau , 2 \tau - \tau^{2} \right)$, a pair of inequalities first conjectured in \cite{HL90}.

The \defn{Kendall distribution function} associated to a general bivariate copula $C$ is defined as the distribution function of the random variable $C(U, V)$, that is,
\[
  K(w) = \proba[C(U, V) \le w], \qquad w \in [0, 1].
\]
The reference to Kendall stems from the link with Kendall's $\tau$, which is given by $\tau = 4 \, \expec[C(U, V)] - 1$. For bivariate Archimedean copulas, for instance, the function $K$ not only identifies the copula \cite{GR93}, convergence of Archimedean copulas is actually equivalent to weak convergence of their Kendall distribution functions \cite{CS08}. For bivariate extreme-value copulas, the function $K$ takes the remarkably simple form
\begin{equation}
\label{E:K}
  K(w) = w - (1 - \tau) \, w \, \log w, \qquad w \in [0, 1],
\end{equation}
as shown in \cite{GKR98}. In fact, in that paper the conjecture was formulated that if the Kendall distribution function of a bivariate copula is given by \eqref{E:K}, then $C$ is a bivariate extreme-value copula, a conjecture which to the best of our knowledge still stands. In the same paper, equation~\eqref{E:K} was used to formulate a test that a copula belongs to the family of extreme-value copulas; see also \cite{BGN09}.

\begin{figure}[h]
\sidecaption[t]
\includegraphics[width=0.5\textwidth]{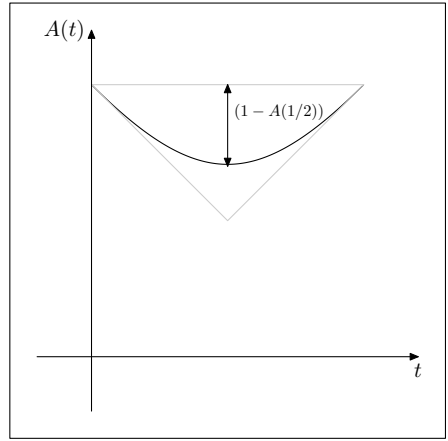} 
\caption{\label{F:lambda} The coefficient of upper tail dependence $\lambda_U$ is equal to twice the length of the double arrow in the upper part of the graph.}
\end{figure}

In the context of extremes, it is natural to study the \defn{coefficient of upper tail dependence}. For a bivariate copula $C_F$ in the domain of attraction of an extreme-value copula with tail dependence function $\ell$ and Pickands dependence function $A$, we find
\begin{align*}
  \lambda_U 
  &= \lim_{u \uparrow 1} \proba(U > u \mid V > u) = \lim_{t \downarrow 0} t^{-1} \bigl( 2t - 1 + C(1-t, 1-t) \bigr) \\
  &= 2 - \ell(1, 1) = 2 \, \bigl(1 - A(1/2) \bigr) \in [0, 1].
\end{align*}
Graphically this quantity can be represented as the length between the upper boundary and the curve of the Pickands dependence function evaluated in the mid-point $1/2$, see Figure~\ref{F:lambda}. The coefficient $\lambda_U$ ranges from $0$ ($A = 1$, independence) to $1$ (complete dependence). Multivariate extensions are proposed in \cite{L09}.

The related quantity $\ell(1, 1) = 2 \, A(1/2)$ is called the \defn{extremal coefficient} in \cite{STY90}. For a bivariate extreme-value copula, we find
\[
  \proba(U \le u, V \le u) = u^{2 \, A(1/2)}, \qquad u \in [0, 1],
\]
so that $2 \, A(1/2) \in [1, 2]$ can be thought of as the (fractional) number of independent components in the copula. Multivariate extensions have been studied in \cite{ST02}.

For the \defn{lower tail dependence coefficient}, the situation is trivial: 
\[
  \lambda_{L} 
  = \lim_{u \downarrow 0 } \proba(U \leq  u \mid V \leq u) =  \lim_{u \downarrow 0 } u^{(2 A(1/2) - 1)} 
  = \begin{cases}
      0 & \text{if $A(1/2) > 1/2$}, \\
      1 & \text{if $A(1/2) = 1/2$}.
     \end{cases}
\]
In words, except for the case of perfect dependence, $A(1/2) = 1/2$, extreme-value copulas have asymptotically independent lower tails.

\section{Estimation}
\label{S:estim}

Let $\mv{X}_i = (X_{i1}, \ldots, X_{id})$, $i \in \{1, \ldots, n\}$, be a random sample from a (continuous) distribution $F$ with margins $F_1, \ldots, F_d$ and extreme-value copula $C$:
\[
  F(x_1, \ldots, x_d) = C\bigl( F_1(x_1), \ldots, F_d(x_d) \bigr),
\]
and $C$ as in Theorem~\ref{T:EVC:H}. The problem considered here is statistical inference on $C$, or equivalently, on its Pickands dependence function $A$. A number of situations may arise, according to whether the extreme-value copula $C$ is completely unknown or is assumed to belong to a parametric family. In addition, the margins may be supposed to be known, parametrically modelled, or completely unknown.

\subsection{Parametric estimation}

Assume that the extreme-value copula $C$ belongs to a parametric family $(C_\theta : \theta \in \Theta)$ with $\Theta \subset \RR^p$; for instance, one of the families described in Section~\ref{S:parametric}. Inference on $C$ then reduces to inference on the parameter vector $\theta$. The usual way to proceed is by maximum likelihood. The likelihood is to be constructed from the copula density
\[
  c_\theta(u_1, \ldots, u_d) = \frac{\partial^d}{\partial u_1 \cdots \partial u_d} C_\theta(u_1, \ldots, u_d), \qquad (u_1, \ldots, u_d) \in (0, 1)^d.
\]
In order for this density to exist and to be continuous, the spectral measure $H$ should be absolute continuous with continuous Radon--Nikodym derivative on all $2^d-1$ faces of the unit simplex with respect to the Hausdorff measure of the appropriate dimension \cite{CT91}. In dimension $d=2$, the Pickands dependence function $A : [0, 1] \to [1/2, 1]$ should be twice continuously differentiable on $(0, 1)$, or equivalently, the spectral measure $H$ should have a continuous density on $(0, 1)$ (after identification of the unit simplex in $\RR^2$ with the unit interval).

In case the margins are unknown, they may be estimated by the (properly rescaled) empirical distribution functions 
\begin{equation}
\label{E:margins}
  \hat{F}_{nj}(x) = \frac{1}{n+1} \sum_{i=1}^n I(X_{ij} \le x), \qquad x \in \RR, \; j \in \{1, \ldots, d\}.
\end{equation}
(The denominator is $n+1$ rather than $n$ in order to avoid boundary effects in the pseudo-loglikelihood below.) Estimation of $\theta$ then proceeds by maximizing the pseudo-loglikelihood
\[
  \sum_{i=1}^n \log c_\theta \bigl( \hat{F}_{n1}(X_{i1}), \ldots, \hat{F}_{nd}(X_{id}) \bigr),
\]
see \cite{GGR95}. The resulting estimator is consistent and asymptotically normal, and its asymptotic variance can be estimated consistently.

If the margins are modelled parametrically as well, a fully parametric model for the joint distribution $F$ arises, and the parameter vector of $F$ may be estimated by ordinary maximum likelihood. An explicit expression for the $5 \times 5$ Fisher information matrix for the bivariate distribution with Weibull margins and Gumbel copula is calculated in \cite{OM1992}. A multivariate extension and with arbitrary generalized extreme value margins is presented in \cite{S95}. 

Special attention to the boundary case of independence is given in \cite{Tawn88}. In this case, the dependence parameter lies on the boundary of the parameter set and the Fisher information matrix is singular, implying the normal assumptions for validity of the likelihood method are no longer valid.

A robustified version of the maximum likelihood estimator is introduced in \cite{DM02}. The effects of misspecification of the dependence structure are studied in \cite{DT01}. 

\subsection{Nonparametric estimation}

For simplicity, we restrict attention here to the bivariate case. For multivariate extensions, see \cite{GdS09, ZWP08}.

Let $(X_1, Y_1), \ldots, (X_n, Y_n)$ be an independent random sample from a bivariate distribution $F$ with extreme-value copula $C$ and Pickands dependence function $A$. Assume for the moment that the marginal distribution functions $F_1$ and $F_2$ are known and put $U_i = F_1(X_i)$ and $V_i = F_2(Y_i)$ and $S_i = - \log U_i$ and $T_i = - \log V_i$. Note that $S_i$ and $T_i$ are standard exponential random variables. For $t \in [0, 1]$, put
\[
  \xi_i(t) = \min \biggl( \frac{S_i}{1 - t}, \frac{T_i}{t} \biggr),
\]
with the obvious conventions for division by zero. A characterizing property of extreme-value copulas is that the distribution of $\xi_i(t)$ is exponential as well, now with mean $1/A(t)$: for $x > 0$,
\begin{align}
\label{E:xit:exp}
  \proba[\xi_i(t) > x]
  &= \proba [U_i < e^{-(1-t)x}, V_i < e^{-tx}] \nonumber \\
  &= C(e^{-(1-t)x}, e^{-tx}) = e^{- x \, A(t)}. 
\end{align}
This fact leads straightforwardly to the original Pickands estimator \cite{Pickands81}:
\begin{equation}
\label{E:Pickands}
  \frac{1}{\hat{A}^{P}(t)} = \frac{1}{n} \sum_{i=1}^n \xi_{i}(t).
\end{equation}
A major drawback of this estimator is that it does not verify any of the constraints imposed on the family of the Pickands dependence functions in Theorem~\ref{T:EVC:A}.

Besides establishing the asymptotic properties of the original Pickands estimator, Deheuvels~\cite{Deheuvels91} proposed an improvement of the Pickands estimator that at least verifies the endpoint constraints $A(0) = A(1) = 1$:
\begin{equation}
\label{E:Deheuvels}
  \frac{1}{\hat{A}^{D}(t)} 
  = \frac{1}{n} \sum_{i=1}^{n} \xi_{i}(t) - t \, \frac{1}{n} \sum_{i=1}^{n} \xi_{i}(1) - (1-t) \, \sum_{i=1}^{n} \xi_{i}(0) + 1.
\end{equation}
As shown in \cite{S07}, the weights $(1-t)$ and $t$ in de Deheuvels estimator \eqref{E:Deheuvels} can be understood as pragmatic choices that could be replaced by suitable weight functions $\beta_{1}(t)$ and $\beta_{2}(t)$:
\begin{equation}
\label{Deheuvels:weight}
  \frac{1}{\hat{A}^{D}(t)} 
  = \frac{1}{n} \sum_{i=1}^{n} \xi_{i}(t) - \beta_{1}(t) \, \frac{1}{n} \sum_{i=1}^{n} \xi_{i}(1) - \beta_{2}(t) \, \sum_{i=1}^{n} \xi_{i}(0) + 1.
\end{equation}
The linearity of the right-hand side of \eqref{E:Deheuvels} in $\xi_i(t)$ suggests to estimate the variance-minimizing weight functions via a linear regression of $\xi_i(t)$ upon $\xi_i(0)$ and $\xi_i(1)$:
\[
  \xi_i(t) = \beta_0(t) + \beta_1(t) \, \{\xi_i(0) - 1\} + \beta_2(t) \, \{ \xi_i(1) - 1 \} + \epsilon_i(t).
\]
The estimated intercept $\hat{\beta}_{0}(t)$ corresponds to the minimum-variance estimator for $1 / A(t)$ in the class of estimators \eqref{Deheuvels:weight}.
 
In the same spirit, Hall and Tajvidi~\cite{HT00} proposed another approach to improve the small-sample properties of the Pickands estimator at the boundary points. For all $t \in [0,1]$ and $i \in \{1, \dots , n \}$, define
\[
  \bar{\xi}_{i}(t) = \min \biggl( \frac{\bar{S}_{i}}{1 - t}, \frac{\bar{T}_{i}}{t} \biggr)
\]
with 
\begin{align*}
  \bar{S}_{i} &= \frac{S_i}{ \frac{1}{n} (S_1 + \cdots + S_n)}, &
  \bar{T}_{i} &= \frac{T_i}{ \frac{1}{n} (T_1 + \cdots + T_n)}.
\end{align*}
The estimator presented in \cite{HT00} is given by
\[
  \frac{1}{\hat{A}^{HT}(t)} = \frac{1}{n} \sum_{i=1}^n \bar{\xi}_{i}(t).
\]
Not only does the estimator's construction guarantee that the endpoint conditions are verified, in addition it always verifies the constraint $\hat{A}^{HT}(t) \ge 1 \vee (1-t)$. Among the three nonparametric estimators mentioned so far, the Hall--Tajvidi estimator typically has the smallest asymptotic variance.

A different starting point was chosen by Cap\'era\`a, Foug\`eres and Genest \cite{CFG97}: they showed that the distribution function of the random variable $Z_i = \log(U_i) / \log(U_i V_i)$ is given by
\[
  \proba (Z_i \le z) = z + z(1-z) \, \frac{A'(z)}{A(z)}, \qquad 0 \le z < 1,
\]
where $A'$ denotes the right-hand derivative of $A$. Solving the resulting differential equation for $A$ and replacing unknown quantities by their sample versions yields the CFG-estimator. In \cite{S07} however, it was shown that the estimator admits the simpler representation
\begin{equation}
\label{E:CFG}
  \log \hat{A}^{CFG}(t) = - \frac{1}{n} \sum_{i=1} \log \xi_i(t) - (1-t) \sum_{i=1}^n \log \xi_i(0) - t \, \sum_{i=1}^n \log \xi_i(1)
\end{equation}
for $t \in [0, 1]$. This expression can be seen as a sample version of
\[
  \expec [- \log \xi_i(t)] = \log A(t) + \gamma, \qquad t \in [0, 1],
\]
a relation which follows from \eqref{E:xit:exp}; note that the Euler--Mascheroni constant $\gamma = 0.5772\ldots$ is equal to the mean of the standard Gumbel distribution. Again, the weights $(1-t)$ and $t$ in \eqref{E:CFG} can be replaced by variance-minimizing weight functions that are to be estimated from the data \cite{GdS09, S07}. The CFG-estimator is consistent and asymptotically normal as well, and simulations indicate that it typically performs better than the Pickands estimator and the variants by Deheuvels and Hall--Tajvidi. \bigskip

Theoretical results for extreme values in the case of \emph{unknown margins} are quite recent. To some extent, Jim\'enez, Villa-Deharce and Flores \cite{RVF01} were the first to present an in-depth treatment of this situation. However, their main theorem on uniform consistency is established under conditions that are unnecessarily restrictive. In \cite{GS09}, asymptotic results were established under much weaker conditions. The estimators are the same as the ones presented above, the only difference being that $U_i = F_1(X_i)$ and $V_i = F_2(Y_i)$ are replaced by
\begin{align*}
  \hat{U}_i &= \hat{F}_{n1}(X_i) = \frac{1}{n+1} \sum_{k=1}^n I(X_k \le X_i), &
  \hat{U}_i &= \hat{F}_{n2}(Y_i) = \frac{1}{n+1} \sum_{k=1}^n I(Y_k \le Y_i),
\end{align*}
with $\hat{F}_{nj}$ as in \eqref{E:margins}. Observe that the resulting estimators are entirely rank-based. Contrary to the case of known margins, the endpoint-corrections are irrelevant in the sense that they do not show up in the asymptotic distribution. Again, the CFG-estimator has the smallest asymptotic variance most of the time. \bigskip

The previous estimators do typically not fulfill the \emph{shape constraints} on $A$ as given in Theorem~\ref{T:EVC:A}. A natural way to enforce these constraints is by modifying a pilot estimate $\hat{A}$ into the convex minorant of $\bigl( \hat{A}(t) \vee (1-t) \vee t \bigr) \wedge 1$, see \cite{Deheuvels91, RVF01, Pickands81}. It can be shown that this transformation cannot cause the $L^\infty$ error of the estimator to increase. A different way to impose the shape constraints is by constrained spline smoothing \cite{AG05, HT00} or by constrained kernel estimation of the derivative of $A$ \cite{STY90}. 

The $L^2$-viewpoint was chosen in \cite{FGS08}. The set $\mathcal{A}$ of Pickands dependence functions being a closed and convex subset of the space $L^2([0, 1], \, \mathrm{d}x)$, it is possible to find for a pilot estimate $\hat{A}$ a Pickands dependence function $A \in \mathcal{A}$ that minimizes the $L^2$-distance $\int_0^1 (\hat{A} - A)^2$. By general properties of orthogonal projections, the $L^2$-error of the projected estimator cannot increase.

Finally, a nonparametric Bayesian approach has been proposed by Guillotte and Perron \cite{GP08}. Driven by a nonparametric likelihood, their methodology yields an estimator with good properties: its estimation error is typically small, it automatically verifies the shape constraints, and it blends naturally with parametric likelihood methods for the margins.

\section{Further reading}
\label{S:further}

About the first monograph to treat multivariate extreme-value dependence is the one by Galambos \cite{Galambos78}, with a major update in the second edition \cite{Galambos87}. Extreme-value copulas are treated extensively in the monographs \cite{BGTS04, NK00} and briefly in the 2006 edition of Nelsen's book \cite{Nelsen06}. The regular-variation approach to multivariate extremes is emphasized in the books by Resnick \cite{Resnick87, Resnick07} and de Haan and Ferreira \cite{dHF06}. A highly readable introduction to extreme-value analysis is the book by Coles~\cite{Coles01}.

The first representations of bivariate extreme-value distributions are due to Finkelstein~\cite{Finkelstein53}, Tiago de Oliveira~\cite{Tiago58}, Geffroy~\cite{Geffroy58, Geffroy59} and Sibuya~\cite{S60}. Incidentally, the 1959 paper by Geffroy appeared in the same issue as the famous paper by Sklar~\cite{Sklar59}. The equivalence of all these representations was shown in Gumbel~\cite{Gumbel62}; see also the more recent paper by Obretenov~\cite{Obretenov91}. However, their representations of multivariate extreme value distributions have not enjoyed the same success as the one proposed by Pickands~\cite{Pickands81}. The domain of attraction condition seems to have been formulated for the first time by Berman~\cite{Berman61}, his standardization being to the standard exponential distribution rather than the uniform one.

A particular class of extreme-value copulas arises if the spectral measure $H$ in Theorem~\ref{T:EVC:H} is discrete. In that case, the stable tail dependence function $\ell$ and the Pickands dependence function $A$ are piecewise linear. In general, such distributions arise from max-linear combinations of independent random variables \cite{ST02}. An early example of such a distribution is the bivariate model studied by Tiago de Oliveira in \cite{TdO74, TdO80, TdO89}, which has a spectral measure with exactly two atoms; see also \cite{EKS08}. The bivariate distribution of Marshall and Olkin~\cite{MO67} corresponds to a spectral measure with exactly three atoms, $\{0, 1/2, 1\}$ (after identification of the unit simplex in $\RR^2$ with the unit interval); see~\cite{MaiScherer09} for a multivariate extension.

Even more challenging than the estimation problem considered in Section~\ref{S:estim} is when the random sample comes from a distribution which is merely in the domain of attraction of a multivariate extreme-value distribution. See for instance \cite{BD07, CT91, dHNP08, EKS08, JSW92, KKP07, KKP08, LT96} for some (semi-)parametric approaches and \cite{AG05, CF00, EdHP01, EdHL06, ES09, SS06} for some nonparametric ones.

For an overview of software related to extreme value analysis, see \cite{SG05}. Particularly useful are the \texttt{R} packages \texttt{evd} \cite{evd}, which provides algorithms for the computation, simulation \cite{Stephenson03} and estimation of certain univariate and multivariate extreme-value distributions, as well as the more general \texttt{copula} package \cite{Yan07}.

\begin{acknowledgement}
The authors' research was supported by IAP research network grant nr.\ P6/03 of the Belgian government (Belgian Science Policy) and by contract nr.\ 07/12/002 of the Projet d'Actions de Recherche Concert\'ees of the Communaut\'e fran\c{c}aise de Belgique, granted by the Acad\'emie universitaire Louvain.
\end{acknowledgement}


\end{document}